\documentclass[11pt,a4paper]{article}
\usepackage{caption2}%---change the caption of figures
\usepackage{CJK}
\usepackage{tabls}
\usepackage{bm}
\usepackage{multirow}
\usepackage{amssymb}
\usepackage{amsfonts}
\usepackage{mathrsfs}
\usepackage{empheq}
\usepackage{amsmath}
\usepackage{amsthm}
\usepackage{graphicx}
\usepackage{color}
\usepackage{indentfirst}
\usepackage{hyperref}%pour PDFLaTex
\usepackage{float}
\usepackage{cases}
\usepackage[titletoc]{appendix}

\allowdisplaybreaks%

 \numberwithin{equation}{section}
\theoremstyle{definition}

 \newtheorem{thm}{Theorem}[section]
 
 \newtheorem{lem}{Lemma}[section]
 \newtheorem{prop}{Proposition}[section]
 \theoremstyle{definition}
 
 \newtheorem{rem}{Remark}[section]

\newcommand{\be}{\begin{equation}}
\newcommand{\ee}{\end{equation}}
\newcommand{\beq}{\begin{equation*}}
\newcommand{\eeq}{\end{equation*}}

\begin{document}

\begin{CJK*}{GB}{gbsn}

\title{\bf  {Boundary feedback stabilization of quasilinear hyperbolic systems with partially dissipative structure}}

\author{ Ke Wang\thanks{Department of Mathematics, Donghua
University, Shanghai 201620, China. E-mail: \texttt{kwang@dhu.edu.cn}. }
\quad   Zhiqiang Wang\thanks{School of Mathematical Sciences and Shanghai Key Laboratory for Contemporary Applied Mathematics,
Fudan University, Shanghai 200433, China.
E-mail: \texttt{wzq@fudan.edu.cn}. }
 \quad Wancong Yao\thanks{School of Mathematical Sciences, Fudan University, Shanghai 200433, China. E-mail: \texttt{wcyao18@fudan.edu.cn}.}
}

\date{March 26, 2020}

\maketitle

\begin{abstract}
In this paper, we study the boundary feedback stabilization of a quasilinear hyperbolic system with partially dissipative structure.
Thanks to this structure, we construct a suitable Lyapunov function which leads to  the exponential stability to the equilibrium of the  $H^2$ solution.
As an application, we also obtain the feedback stabilization for the Saint-Venant-Exner model under physical boundary conditions.
 \end{abstract}



\section{\hspace{-5mm}.\ Introduction and main results}

\noindent

Many models  in physics, mechanics and other fields, including  gas and fluid dynamics for instance, are described as hyperbolic equations. 
Control problems, particularly the stability and stabilization problems,  of hyperbolic systems have  been widely studied for decades 
(see \cite{Coron2016book,CoronBook07} and the references therein).

One classical approach to establish the asymptotic stability of hyperbolic system is the characteristic method. In the framework of $C^1$-solution, 
dissipative boundary conditions that lead to exponential stability  of quasilinear  hyperbolic systems without source terms have been found in \cite{Greenberg1984effect, Li1994book}.

Another important approach to design boundary feedback controls is  the Backstepping method.  It has been used to stabilize exponentially  the inhomogeneous 
quasilinear hyperbolic system in $H^2$ norm (see \cite{Coron2013backstepping,HRDK19}).  
One can refer to  \cite{Krstic2008book} for many successful examples about feedback stabilization  with this approach.

The third powerful approach is the Lyapunov function method.  A strict Lyapunov function is introduced in \cite{Xu2002} to 
achieve the exponential stability of a class of symmetric linear hyperbolic systems. Similar  Lyapunov functions are used for 
quasilinear homogeneous hyperbolic systems in the framework of $H^2$-solution in \cite{Coron2008dissipative}.

If the hyperbolic system is inhomogeneous,  the  Lyapunov function approach can still be applied (see for instance \cite{Dick2011strict,Gugat2017}).
However, the nonzero source term change a lot the stability properties.
With a source term,    a simple quadratic Lyapunov function ensuring exponential stability for the $L^2$ norm (or $H^p$ norm) does not always exist
 no matter what the boundary conditions are.
In \cite{BC11}, the authors study a linear $2\times2$ hyperbolic system and found a necessary and sufficient condition for simple quadratic
Lyapunov function.  Later in Chapter 6 of \cite{Coron2016book},  the authors give a sufficient (but a priori non-necessary) condition such that 
the exponential stability of the system for the $H^p$ norm with $p\geq 2$ is achieved.  
We refer to \cite{Hayat2018} for a relevant result in $C^1$-norm or $C^p$-norm.
Naturally, these conditions all include one interior condition which requires \emph{a good coupling structure} of the hyperbolic system, compared to the homogeneous case.
However,  as mentioned in their papers,  this interior condition (typically a differential matrix inequality) is not straightforward to be checked in a specific model.

%
%[2] Georges Bastin and Jean-Michel Coron. On boundary feedback stabilization of non-uniform linear 2
%2 hyperbolic systems over a bounded interval. Systems & Control Letters, 60(11):900{906, 2011.

%

%

Different from the above, 
Herty and Yong study the boundary feedback stabilization of one-dimensional \emph{linear} hyperbolic systems with a relaxation term in  \cite{Yong2016}.
The key assumption is a \emph{structural stability condition} which is introduced from \cite{Yong1999} and is satisfies in many physical models.
Later, in \cite{Yong2019}, Yong shows that under this \emph{structural stability condition}, 
the boundary feedback stabilization result is also available for a class of  one-dimensional linear hyperbolic system with vanishing eigenvalues.

Motivated by \cite{Yong2016}, in this paper, we consider a one-dimensional quasilinear hyperbolic system with the same relaxation structure. 
Thanks to the partial dissipation in the  \emph{structural stability condition}, we  establish the local exponential stability of this \emph{nonlinear} system for the $H^2$-norm. 
The main strategy is to construct a strict  Lyapunov function together with a perturbation argument based on linear approximation.
Compared to  the result in \cite{Yong2016}, we provide an explicit sufficient condition on the gains of stabilizing boundary feedback control.
As an application, we also obtain the boundary feedback stabilization of the Saint-Venant-Exner model proposed in \cite{Hudson2003} under physical boundary conditions.

Precisely, we are concerned with the boundary feedback stabilization of the following one-dimensional quasilinear hyperbolic system
\begin{align}\label{system}
U_t + \mathbf{A}(U)U_x = \mathbf{Q}(U), \quad t \in (0,\infty), x\in (0,1)
\end{align}
where $U=(u_1,\cdots,u_n)^T$ is the unknown vector function of $(t,x)$,
$\mathbf{A}: \mathbb{R}^n \mapsto \mathcal{M}_{n,n}(\mathbb{R})$ is a smooth matrix function and  $\mathbf{Q}: \mathbb{R}^n \mapsto \mathbb{R}^{n}$ is a smooth vector function.

Let $U^*\in \mathbb{R}^n$ be an equilibrium of \eqref{system}, i.e.,
\begin{align}\label{steady}
\mathbf{Q}(U^*) = 0.
\end{align}
Without loss of generality, we may assume  $U^*=0$, otherwise one can consider $U-U^*$ as the unknown functions.

We first assume that the system \eqref{system} is hyperbolic in a neighborhood of $U=0$, i.e., the matrix $\mathbf{A}(U)$ has $n$ real eigenvalues
\begin{align}\label{eigenvalues}
\mathbf{\Lambda}_r(U)<0<\mathbf{\Lambda}_s(U)\quad (r=1,\cdots,m;\ s=m+1,\cdots,n),
\end{align}
and it has a complete set of left eigenvectors $\mathbf{L}_i(U)=(\mathbf{L}_{i1}(U),\cdots, \mathbf{L}_{in}(U)), \ (i=1,\cdots,n)$, i.e.,
\begin{align}\label{eigenvec}
\mathbf{L}_i(U)\mathbf{A}(U) = \mathbf{\Lambda}_i(U) \mathbf{L}_i(U)    \quad (i=1,\cdots,n).
\end{align}
Let
\begin{align}\label{matrixtL}
\mathbf{L}(U)=\left(\begin{array}{ll}\mathbf{L}_1(U)\\ \vdots \\ \mathbf{L}_n(U)
\end{array}\right)
\text{\quad and \quad }
\mathbf{\Lambda}(U)=\left(\begin{array}{cc}
     \mathbf{\Lambda}_-(U) & 0\\
     0 & \mathbf{\Lambda}_+(U)
\end{array}\right)
\end{align}
where
\begin{align}
\mathbf{\Lambda}_-(U) =\text{diag}\{\mathbf{\Lambda}_1(U),\cdots,\mathbf{\Lambda}_m(U)\},
\text{\quad and \quad }
\mathbf{\Lambda}_+(U) =\text{diag}\{\mathbf{\Lambda}_{m+1}(U),\cdots,\mathbf{\Lambda}_n(U)\}.
\end{align}
Then
\begin{align}\label{LAL-1}
\mathbf{L}(U)\mathbf{A}(U) = \mathbf{\Lambda}(U) \mathbf{L}(U).
\end{align}

It is easy to see that system \eqref{system} is hyperbolic if and only if
there is a symmetric positive definite matrix $\mathbf{A}_0(U)$, such that
\begin{align}\label{hyperbolic}
\mathbf{A}_0(U)\mathbf{A}(U) = \mathbf{A}^T(U)\mathbf{A}_0(U).
\end{align}
Then it follows that
\begin{align}
(\mathbf{L}^{-1}(0))^T\mathbf{A}_0(0)\mathbf{L}^{-1}(0)\mathbf{\Lambda}(0)
= \mathbf{\Lambda}(0) (\mathbf{L}^{-1}(0))^T\mathbf{A}_0(0)\mathbf{L}^{-1}(0).
\end{align}
Consequently, there exist two symmetric positive definite matrices
$\mathbf{X}_1(0)\in \mathcal{M}_{m,m}(\mathbb{R})$ and $\mathbf{X}_2(0)\in \mathcal{M}_{n-m,n-m}(\mathbb{R})$ such that
\begin{align}\label{block-diagonal}
(\mathbf{L}^{-1}(0))^T\mathbf{A}_0(0)\mathbf{L}^{-1}(0) = \left(\begin{array}{cc}
     \mathbf{X}_1(0) & 0\\
     0 & \mathbf{X}_2(0)
\end{array}\right).
\end{align}

Moreover,  we assume the system possesses the following \emph{partially dissipative structure} in a neighborhood of $U=0$:

 There exist invertible matrices $\mathbf{P}(U)\in \mathcal{M}_{n,n}(\mathbb{R})$ and $\mathbf{S}(U)\in \mathcal{M}_{r,r}(\mathbb{R})$  with  $0<r  \leq n$, such that
\begin{align}\label{structural1}
&\mathbf{P}(U)\mathbf{Q}_U(U) = \left(\begin{array}{cc}
     0 & 0\\
     0 & \mathbf{S}(U)
\end{array}\right)\mathbf{P}(U),
 	\\
\label{structural2}
&\mathbf{A}_0(U)\mathbf{Q}_U(U) + \mathbf{Q}_U^T(U)\mathbf{A}_0(U) \leq -\mathbf{P}^T(U)
\left(\begin{array}{cc}
     0 & 0\\
     0 & \mathbf{I}_r
\end{array}\right)\mathbf{P}(U).
\end{align}
Here $\mathbf{Q}_U(U)$ stands for the Jacobian matrix of $\mathbf{Q}$ with respect to $U$,  $I_r$ denotes the $r\times r$ identity matrix.
Let us point out  that the above assumptions \eqref{structural1} and \eqref{structural2}  are called \emph{structural stability conditions}  in \cite{Yong1999, Yong2016},
which are commonly satisfied in lots of physical models.

Let
\begin{align}\label{bcvariable}
\xi(t, x) = \mathbf{L}(0)U(t, x)
\end{align}
be  the  linearized diagonal variable and denote
\begin{align} \label{xi+-}
\xi(t,x)=
\left (\begin{array}{ll}
\xi_-(t, x)\\
\xi_+(t, x)
\end{array}\right)
\end{align}
where $\xi_-(t,x)\in \mathbb{R}^m$ and $\xi_+(t,x)\in \mathbb{R}^{n-m}$.
According to the theory on the well-posedness of the quasilinear hyperbolic system,
the typical boundary conditions are given as follows
\begin{align}\label{bc}
\left(\begin{array}{ll}
\xi_+(t, 0)\\
\xi_-(t, 1)
\end{array}\right)
= \mathbf{K}\left(\begin{array}{ll}
\xi_+(t, 1)\\
\xi_-(t, 0)
\end{array}\right), \quad  t\in (0,\infty),
\end{align}
where the feedback matrix
\begin{align}\label{matrixK}
 \mathbf{K} = \left(\begin{array}{ll}
     \mathbf{K}_{00} & \mathbf{K}_{01}\\
     \mathbf{K}_{10} & \mathbf{K}_{11}
\end{array}\right),
\end{align}
 $\mathbf{K}_{00}\in \mathcal{M}_{n-m,n-m}(\mathbb{R})$, $\mathbf{K}_{01}\in \mathcal{M}_{n-m,m}(\mathbb{R})$, $\mathbf{K}_{10}\in \mathcal{M}_{m,n-m}(\mathbb{R})$ and $\mathbf{K}_{11}\in \mathcal{M}_{m,m}(\mathbb{R})$ are all matrices with constant elements.

Finally, the initial condition is prescribed as
\begin{align}\label{id}
U(0, x) = U_0(x),\quad x\in (0,1),
\end{align}
with $U_0 \in H^2((0, 1);\mathbb{R}^n)$  in a neighborhood of $U=0$.

Regarding the well-posedness of the solutions to the problem \eqref{system}, \eqref{bc} and \eqref{id}, we
have the following proposition
\begin{prop}\label{prop}
There exists $\delta_0 > 0$ such that, for every $U_0 \in H^2((0,1);\mathbb{R}^n)$
satisfying
\begin{align}
||U_0||_{H^2((0,1);\mathbb{R}^n)}\leq \delta_0,
\end{align}
and the $C^1$ compatibility conditions at the points $(t,x)=(0,0)$, $(0,1)$, the problem \eqref{system}, \eqref{bc} and \eqref{id} has a unique maximal classical solution
\begin{align}
U\in C^0([0,T);H^2((0,1);\mathbb{R}^n))
\end{align}
with $T\in(0,+\infty]$. Moreover, if
\begin{align}
||U(t,\cdot)||_{H^2((0,1);\mathbb{R}^n)}\leq \delta_0,\quad \forall t\in[0,T),
\end{align}
then $T=+\infty$.
\end{prop}

\vspace{2mm}
Our main result is the following theorem.

\begin{thm}\label{thm}
Assume that the hyperbolic system \eqref{system} has no vanishing characteristic speed and it possesses the partially dissipative structure, i.e.,
 \eqref{eigenvalues}, \eqref{hyperbolic}, \eqref{structural1} and \eqref{structural2} hold.
Let  $\mathbf{K}$ be chosen such that matrices
\begin{align}\label{K1}
\left(\begin{array}{cc}
      \mathbf{X}_2(0)\mathbf{\Lambda}_+(0) & 0\\
      0 &- \mathbf{X}_1(0)\mathbf{\Lambda}_-(0)
\end{array}\right)- \mathbf{K}^T\left(\begin{array}{cc}
     \mathbf{X}_2(0)\mathbf{\Lambda}_+(0) & 0\\
      0 & -\mathbf{X}_1(0)\mathbf{\Lambda}_-(0)
\end{array}\right)\mathbf{K}
\end{align}
and
\begin{align}\label{K2}
\left(\begin{array}{cc}
     e^{-\mathbf{\Lambda}_+(0)}\mathbf{\Lambda}_+(0) & 0\\
      0 & -\mathbf{\Lambda}_-(0)
\end{array}\right)- \mathbf{K}^T\left(\begin{array}{cc}
     \mathbf{\Lambda}_+(0) & 0\\
      0 & -e^{-\mathbf{\Lambda}_-(0)}\mathbf{\Lambda}_-(0)
\end{array}\right)\mathbf{K}
\end{align}
are both positive definite.
Then,  the closed-loop system
  \eqref{system}, \eqref{bc} and \eqref{id} is locally exponentially stable for the $H^2$-norm, i.e., there exist positive constants $\delta$, $C$ and $\nu$, such that
the solution to the system \eqref{system}, \eqref{bc} and \eqref{id} satisfies
\begin{align}\label{exponentially stable}
||U(t,\cdot)||_{H^2((0, 1);\mathbb{R}^n)} \leq Ce^{-\nu t}||U_0||_{H^2((0, 1);\mathbb{R}^n)},\quad  t \in [0, +\infty),
\end{align}
provided that
\begin{equation}\label{eU0}
||U_0||_{H^2((0, 1);\mathbb{R}^n)} \leq \delta
\end{equation}
and the $C^1$ compatibility conditions are satisfied at $(t,x)=(0,0)$ and $(0,1)$.
\end{thm}

\begin{rem}\label{rem2}
Theorem \ref{thm} still holds if  the assumption \eqref{structural1} is extended to a more general case that
\begin{align}\label{structural12}
\mathbf{P}(U)\mathbf{Q}_U(U)\mathbf{P}^{-1}(U)=\left(\begin{array}{cc}
\mathbf{S}_{11}(U) & \mathbf{S}_{12}(U)\\
\mathbf{S}_{21}(U) & \mathbf{S}_{22}(U)
\end{array}\right),
\end{align}
where $\mathbf{S}_{22}(U)\in \mathcal{M}_{r,r}(\mathbb{R})$ with $0 < r \leq n$ is an invertible matrix,
$|\mathbf{S}_{11}(0)|_{\infty}$ and $|\mathbf{S}_{21}(0)|_{\infty}$ are sufficiently small.
\end{rem}

\begin{rem}\label{rem3}
Theorem \ref{thm} still holds if the assumption \eqref{structural2} is extended to a more general case that
\begin{align}\label{structural22}
\mathbf{A}_0(U)\mathbf{Q}_U(U) + \mathbf{Q}_U^T(U)\mathbf{A}_0(U) \leq -\mathbf{P}^T(U)
\left(\begin{array}{cc}
     0 & 0\\
     0 & \mathbf{R}(U)
\end{array}\right)\mathbf{P}(U),
\end{align}
where $\mathbf{R}(U)\in\mathcal{M}_{r,r}(\mathbb{R})$ is a symmetric positive definite matrix.
\end{rem}

\begin{rem}\label{rem5}
The conditions on the feedback matrix  \eqref{K1} and \eqref{K2} are satisfied provided that
 $|\mathbf{K}|_{\infty}$ is  sufficiently small.
Particularly,  if $\mathbf{K}$ is chosen as \eqref{matrixK} where
\begin{align*}
\mathbf{K}_{00} = \kappa_+ \mathbf{I}_{n-m},  \quad \mathbf{K}_{11} = \kappa_- \mathbf{I}_{m}, \quad \mathbf{K}_{01} = 0, \quad \mathbf{K}_{10} = 0,
\end{align*}
with two  constants $\kappa_{+}$ and $\kappa_{-}$ satisfying
\begin{align*}
\kappa_{+}^2 < \exp\Big(-\max\limits_{s=m+1,\cdots,n}\mathbf{\Lambda}_s(0)\Big)\quad \text{and}\quad
\kappa_-^2 < \exp\Big(\min\limits_{r=1,\cdots,m} \mathbf{\Lambda}_r(0)\Big).
\end{align*}
\end{rem}

\begin{rem}
Let us emphasize that the partially dissipative structure \eqref{structural1} and \eqref{structural2} combined with the dissipative boundary conditions \eqref{K1} and \eqref{K2}
can be included by the interior and boundary stability conditions proposed in \cite[Theorem 6.6]{Coron2016book}.  
However,  the interior conditions on stability are typically differential matrix inequality, 
while the conditions proposed in this paper are all algebraic conditions  which are more straightforward to be checked. 
\end{rem}

\begin{rem}
It is also worthy of mentioning that  the stability conditions (both interior and boundary conditions) for the nonlinear hyperbolic systems 
depend on the topology and  in particular that the stability in $H^2$ norm does not imply the stability in $C^1$ norm (see \cite{CN15}).
\end{rem}

The paper is organized as follows: in Section 2, we introduce  a new quasilinear hyperbolic system with a simpler structure under a  transformation of the unknown functions. Then in Section 3, we construct  a weighted $H^2$-Lyapunov function to prove the exponential stability of the new system which implies immediately Theorem \ref{thm}.  The proofs of related lemmas will be given in Section 4. Finally, in Section 5, the main result is applied to the  Saint-Venant-Exner model for moving water on an open canal under physical boundary conditions.

\section{Transformation of the system}

In this section, we introduce  a new hyperbolic system with a \emph{partially dissipative but simpler structure} under a  transformation of the unknown functions.
In this way, the exponential stability of the original system is reduced to that of the new system.

Let
\begin{align}\label{variableV}
V = \mathbf{P}(0)U.
\end{align}
Then, the  system \eqref{system} can be reduced to
\begin{align}\label{systemV}
V_t + A(V)V_x = B(V),
\end{align}
where
\begin{align}\label{matrixAB}
A(V) = \mathbf{P}(0)\mathbf{A}(\mathbf{P}^{-1}(0)V)\mathbf{P}^{-1}(0) \text{\quad and \quad} B(V)=\mathbf{P}(0)\mathbf{Q}(\mathbf{P}^{-1}(0)V).
\end{align}
Clearly, $V=0$ is an equilibrium of \eqref{systemV} and the Jacobian matrix of $B$ with respect to $V$ at the equilibrium can be calculated as
\begin{equation} \label{BV}
B_V(0)=\mathbf{P}(0)\mathbf{Q}_U(0)\mathbf{P}^{-1}(0).
\end{equation}

Let
\begin{align}\label{matrixL}
L(V) = \mathbf{L}(\mathbf{P}^{-1}(0)V)\mathbf{P}^{-1}(0)
 \text{\quad and \quad}
\Lambda(V) = \mathbf{\Lambda}(\mathbf{P}^{-1}(0)V).
\end{align}
It is easy to check that $L(V)$ is the matrix composed  of the left eigenvectors of $A(V)$, i.e.,
\begin{align}\label{matrixA2}
L(V) A(V) = \Lambda(V)L(V),
\end{align}
which implies that system \eqref{systemV} is a hyperbolic system without vanishing characteristic speeds. Obviously, we have
\begin{align}
\Lambda(0)=\mathbf{\Lambda}(0)
=\text{diag}\{\mathbf{\Lambda}_1(0),\cdots,\mathbf{\Lambda}_n(0)\}.
\end{align}

Let
\begin{align}\label{A0V}
A_0(V)=(\mathbf{P}^{-1}(0))^T\mathbf{A}_0(\mathbf{P}^{-1}(0)V)\mathbf{P}^{-1}(0).
\end{align}
Obviously, $A_0(V)$ is a symmetric positive definite matrix satisfying
\begin{align}\label{A0A}
A_0(V)A(V) = A^T(V)A_0(V).
\end{align}

%and in particular
%can be rewritten as
%\begin{align}\label{matrixB}
%B(V) = B_V(0)V+ \mathcal{O}(V^2),
%\end{align}
%
%
%Here and hereafter $\mathcal{O}(V^2)$ denotes some term of order $V^2$.

Thanks to \eqref{BV} and \eqref{A0V},  the \emph{partially dissipative structure} \eqref{structural1}-\eqref{structural2}  for the original system \eqref{system}
implies  the following \emph{partially dissipative but simpler structure} for system  \eqref{systemV} at  the equilibrium $V=0$.
\begin{align}\label{matrixBV}
& B_V(0)=\left(\begin{array}{cc}
     0 & 0\\
     0 & \mathbf{S}(0)
\end{array}\right),
\\\label{A0B}
&A_0(0)B_V(0) + B^T_V(0)A_0(0) \leq -\left(\begin{array}{cc}
     0 & 0\\
     0 & \mathbf{I}_r
\end{array}\right).
\end{align}
According to the structure \eqref{matrixBV} and \eqref{A0B}, we write $V(t,x)$  as
\begin{align} \label{v1v2}
V(t,x) = \left(\begin{array}{c}
v_1(t,x)\\
v_2(t,x)
\end{array}\right)
\quad  \text {with \ }v_1\in \mathbb{R}^{n-r},\  v_2\in \mathbb{R}^r
\end{align}
for further use.

From \eqref{bcvariable} and \eqref{matrixL}, we can see that the linear diagonal variables $\xi(t,x)$ now becomes
\begin{align}\label{bcvariable1}
\xi(t,x) =  L(0)V(t,x),
\end{align}
with \eqref{xi+-},  which implies that the boundary conditions are still given by \eqref{bc}.

The initial condition for the variable $V$ is given by
\begin{align}\label{Vid}
V(0, x) = V_0(x)\triangleq\mathbf{P}(0)U_0(x), \quad x\in (0,1).
\end{align}

In order to prove Thereom \ref{thm}, it suffices to establish the $H^2$-stabilization for the system  \eqref{systemV}, \eqref{Vid} and  \eqref{bc}.

\section{Proof of Theorem \ref{thm}}

In this section, we can find suitable conditions on the feedback matrix $\mathbf{K}$ such that  the closed-loop system \eqref{systemV}, \eqref{Vid} and  \eqref{bc} is exponentially stable in $H^2(0,1)$-norm.  Then Theorem \ref{thm} follows immediately.

Let $V_0$ with small $H^2((0,1);\mathbb{R}^n)$ norm be such that  the $C^1$ compatibility conditions at $(t,x)=(0,0)$ and $(0,1)$ are satisfied.
Let also  $V\in C^0([0,T),H^2((0,1);\mathbb{R}^n))$ be the maximal classical solution of the
problem \eqref{systemV}, \eqref{Vid} and \eqref{bc}.
We remark that we only prove the stabilization result for smooth solutions while the conclusion follows easily from an density and continuity arguments for distributed solutions.

Motivated by \cite{Coron2016book} and \cite{Yong2016}, we construct a weighted Lyapunov function as follows:
\begin{align}\label{LF}
\mathbb{L}(t) \triangleq  \mathbb{L}_0(t) + \mathbb{L}_1(t) + \mathbb{L}_2(t)
\end{align}
with
\begin{eqnarray}\label{L0}
\mathbb{L}_0(t) \triangleq  \int^1_0V^T\Big(\alpha A_0(V) + L^T(V)e^{-\mathbf{\Lambda}(0)x}L(V)\Big)V\,\mathrm{d}x,\\ \label{L1}
\mathbb{L}_1(t) \triangleq  \int^1_0V_t^T\Big(\alpha A_0(V) + L^T(V)e^{-\mathbf{\Lambda}(0)x}L(V)\Big)V_t\,\mathrm{d}x,\\ \label{L2}
\mathbb{L}_2(t) \triangleq  \int^1_0V_{tt}^T\Big(\alpha A_0(V) + L^T(V)e^{-\mathbf{\Lambda(0)}x}L(V)\Big)V_{tt}\,\mathrm{d}x,
\end{eqnarray}
where  $e^{-\mathbf{\Lambda}(0)x}=\text{diag}\{e^{-\mathbf{\Lambda}_1(0)x},\cdots,e^{-\mathbf{\Lambda}_n(0)x}\}$, $\alpha>0$ is a constant to be chosen later.

For the simplicity of statements, we denote the $\|\cdot\|_{L^2(0,1)}$ norm as $\|\cdot\|$,
$\|\cdot\|_{C^0([0,1])}$ norm as $|\cdot|_0$, $\|\cdot\|_{C^1([0,1])}$ norm as $|\cdot|_1$.

By definition of the Lyapunov function $\mathbb{L}(t)$, $\mathbb{L}(t)$ is equivalent to the energy $||V||^2+||V_t||^2+||V_{tt}||^2$  if $|V(t,\cdot)|_0$ is small. On the other hand,
Differentiation of system \eqref{systemV} with respect to $t$ and $x$ gives that
\begin{align}\label{Vtt}
&V_{tt} + A(V)V_{tx} = B_V(V)V_t - (A'(V)V_t)V_x,\\ \label{Vtx}
&V_{tx} + A(V)V_{xx} = B_V(V)V_x - (A'(V)V_x)V_x,
\end{align}
in which $A'(V)V_t, A'(V)V_x$ are matrices with entries $\frac{\partial a_{ij}(V)}{\partial V}V_t, \frac{\partial a_{ij}(V)}{\partial V}V_x$ respectively.
Then   it is easy to see that if $\mathbb{L}(t)$ is equivalent to the energy $||V(t, \cdot)||^2_{H^2((0,1);\mathbb{R}^n)}$ if $|V(t,\cdot)|_1$ is small.

Next, we turn to estimate  the time derivative of $\mathbb{L}(t)$. For this purpose, we can establish the following  lemmas
with the assumptions \eqref{K1}-\eqref{K2}  on the feedback matrix $\mathbf{K}$.  The proof of the lemmas will be given in Section 4.
\begin{lem}\label{lem1}
There exist positive constants $\alpha_0$, $\beta_0$, $\gamma_0$ and $\delta_0$ independent of $V$ such that if $|V(t,\cdot)|_0\leq \delta_0$,
\begin{align}\label{el0}
\mathbb{L}_0'(t) \leq - \alpha_0 ||v_1||^2 + (\beta_0 - \alpha)||v_2||^2  + \gamma_0 \int^1_0 (|V|^3 + |V|^2|V_t|)\,\mathrm{d}x,
\end{align}
where $v_1,v_2$ is defined in \eqref{v1v2}.
\end{lem}

\begin{lem}\label{lem2}
There exist positive constants $\alpha_1$, $\beta_1$, $\gamma_1$ and $\delta_1$ independent of $V$ such that if $|V(t,\cdot)|_0\leq \delta_1$,
\begin{align}\label{el1}
\mathbb{L}_1'(t) \leq - \alpha_1||v_{1t}||^2 + (\beta_1 - \alpha)||v_{2t}||^2 +  \gamma_1 \int^1_0 (|V_t|^3+|V||V_t|^2) \,\mathrm{d}x.
\end{align}
\end{lem}

\begin{lem}\label{lem3}
There exist positive constants $\alpha_2$, $\beta_2$, $\gamma_2$ and $\delta_2$ independent of $V$ such that if $|V(t,\cdot)|_1\leq \delta_2$,
\begin{align}\label{el2}
\mathbb{L}_2'(t) \leq - \alpha_2||v_{1tt}||^2 + (\beta_2 - \alpha)||v_{2tt}||^2
+  \gamma_2 \int^1_0 (|V_t||V_{tt}|^2 + |V||V_{tt}|^2 + |V_t|^2 |V_{tt}|)\,\mathrm{d}x.
\end{align}
\end{lem}

\vspace{2mm}

With the help of Lemmas \ref{lem1}, \ref{lem2} and \ref{lem3}, we are ready to prove Theorem \ref{thm}.

Let the constants $\alpha> \max\{\beta_0,\beta_1,\beta_2\}$ and $\delta_4  \leq \min\{\delta_0,\delta_1,\delta_2,\delta_3\}$.
The combination of \eqref{el0}-\eqref{el2} yields that there exist positive constants $\beta$ and $\gamma$ such that
\begin{align}
\mathbb{L}'(t)& \leq -\beta \mathbb{L}(t)+ \gamma |V(t,\cdot)|_1\mathbb{L}(t),
\end{align}
if $|V(t,\cdot)|_1\leq \delta_4$.
%\begin{align} \label{H}
%&\mathbb{H} =
%\int^1_0  \big(|V| + |V_t|\big) \big(|V|^2 +|V_t|^2 + |V_{tt}|^2 \big)\mathrm{d}x,
%\end{align}
%
%

Let $\delta_5 \triangleq \min\{\delta_4, \frac{\beta}{2\gamma}\}$.
If we assume \emph{in a priori} that $|V(t,\cdot)|_{1} \leq  \delta_5$ for $t\in (0,T)$,  we get
\begin{align}\label{L't}
\mathbb{L}'(t) \leq -\frac{\beta}{2}\mathbb{L}(t), \quad t\in (0,T)
\end{align}
which implies that $\mathbb{L}(t)$ decays exponentially
\begin{align}
\mathbb{L}(t) \leq  e^{-\frac{\beta t}{2}}\mathbb{L}(0),  \quad t\in (0,T).
\end{align}
Using the equivalence of  the energy $||V(t, \cdot)||^2_{H^2((0,1);\mathbb{R}^n)}$ and $ \mathbb{L}(t)$, we obtain
\begin{align}\label{Vdecay}
||V(t, \cdot)||_{H^2((0,1);\mathbb{R}^n)} \leq C_{1} e^{-\frac{\beta t}{2}}||V_0||_{H^2((0,1);\mathbb{R}^n)},\quad \forall t\in [0,T)
\end{align}
for some constant $C_{1}>0$.

Note also the Sobolev inequality implies
\begin{align}\label{V1H2}
|V(t,\cdot)|_1 \leq C_2||V(t,\cdot)||_{H^2((0,1);\mathbb{R}^n)} \leq C_2C_{1} ||V_0||_{H^2((0,1);\mathbb{R}^n)},\quad \forall t\in [0,T).
\end{align}

Let now
\begin{equation}
\delta= \frac{\delta_5}{C_2 C_1}.
\end{equation}
then the a priori estimate on $|V(t,\cdot)|_1 \leq \delta_5$ indeed holds in $[0,T)$ if  $||V_0||_{H^2((0,1);\mathbb{R}^n)}\leq \delta$.   Therefore
\eqref{Vdecay} follows immediately.
According to Proposition \ref{prop}, we finally conclude that the inequality \eqref{Vdecay} indeed holds for $T=+\infty$.
%\begin{align}\label{Vdecay}
%||V(t, \cdot)||_{H^2((0,1);\mathbb{R}^n)} \leq C_2^2 e^{-\frac{\beta t}{2}}||V_0||_{H^2((0,1);\mathbb{R}^n)},\quad \forall t\geq0.
%\end{align}

Consequently, it follows from \eqref{variableV} that the solution $U$ to the problem \eqref{system}, \eqref{bc} and \eqref{id} is locally  exponentially stable for the $H^2$-norm.
This concludes the proof of Theorem \ref{thm}.

\section{Proof of the Lemmas}

\subsection{Proof of Lemma 2.1}
We calculate the time-derivative of $\mathbb{L}_0(t)$ defined by \eqref{L0},
\begin{align}\label{L0'}
\mathbb{L}_0'(t)=\int ^1_0 2\Big( \alpha V^TA_0(V)V_t +  V^T L^T(V)e^{-\mathbf{\Lambda}(0)x}L(V)V_t\Big)\,\mathrm{d}x
+ \mathcal{O}\Big(\int^1_0 |V|^2|V_t|\,\mathrm{d}x; |V(t,\cdot)|_0\Big).
\end{align}
Here  and hereafter $\mathcal{O}(X;Y)$  denotes the terms that for $X\geq 0$, $Y\geq0$, there exist $C>0$ and $\varepsilon>0$ independent of $V$, $V_t$ and $V_{tt}$, satisfying
\begin{align}
Y\leq \varepsilon \Rightarrow |\mathcal{O}(X;Y)| \leq CX.
\end{align}

Substituting the system \eqref{systemV} into \eqref{L0'}, we have
\begin{align}\label{L0'2}
\mathbb{L}_0'(t) =  \mathbb{I}_0+ \mathbb{J}_0+ \mathcal{O}\Big(\int^1_0 |V|^2|V_t|\,\mathrm{d}x; |V(t,\cdot)|_0\Big)
\end{align}
where
\begin{align}\label{I}
&  \mathbb{I}_0 \triangleq \int^1_0  2\alpha \Big(V^TA_0(V)B(V) - V^TA_0(V)A(V)V_x \Big) \,\mathrm{d}x,\\\label{II}
&  \mathbb{J}_0\triangleq \int^1_0  2\Big(V^T L^T(V)e^{-\mathbf{\Lambda}(0)x}L(V)B(V) -  V^T L^T(V)e^{-\mathbf{\Lambda}(0)x}L(V)A(V)V_x\Big)\,\mathrm{d}x.
\end{align}
Let's first estimate  the term $ \mathbb{I}_0$.
Using \eqref{A0A}, \eqref{A0B} and  integrations by parts, we have
\begin{align} \nonumber
  \mathbb{I}_0 &=\int^1_0 \alpha V^T\Big(A_0(0)B_V(0) + B^T_V(0)A_0(0) \Big)V - \Big(\alpha V^TA_0(V)A(V)V\Big)_x\,\mathrm{d}x
\\ \nonumber &\quad +\mathcal{O}\Big(\int^1_0 (|V|^3+|V|^2|V_t|) \,\mathrm{d}x; |V(t,\cdot)|_0\Big)
\\  \label{I1}
&\leq -\alpha||v_2||^2 + \Big[-\alpha V^T A_0(V)A(V)V\Big]\Big|_0^1 +\mathcal{O}\Big(\int^1_0 (|V|^3+|V|^2|V_t|) \,\mathrm{d}x; |V(t,\cdot)|_0\Big).
\end{align}

%Integration by parts for the term $\Big(V^TA^T(V)A_0(V)V\Big)_x$ yields
%\begin{align}\label{I}
%I\leq  -||v_2||^2  - V^T A_0(0)A(0)V|^1_0 + \int^1_0  O(V^3 + V^2V_t) \mathrm{d}x
%\end{align}

Then we turn to estimate the term $ \mathbb{J}_0$. Linear approximation together with \eqref{matrixA2} implies
\begin{align} \nonumber
 \mathbb{J}_0&= \int^1_0 2V^T L^T(0)e^{-\mathbf{\Lambda}(0)x}L(0)B_V(0)V\,\mathrm{d}x
- \int^1_0 V^TL^T(0)e^{-\mathbf{\Lambda}(0)x}\mathbf{\Lambda}^2(0)L(0)V\,\mathrm{d}x  \\ \label{J0}
&\quad -  \int^1_0\Big(V^TL^T(V)e^{-\mathbf{\Lambda}(0)x}\Lambda(V)L(V)V\Big)_x\,\mathrm{d}x +\mathcal{O}\Big(\int^1_0 (|V|^3+|V|^2|V_t|) \,\mathrm{d}x; |V(t,\cdot)|_0\Big).
\end{align}

Note that  the matrix $L^T(0)e^{-\mathbf{\Lambda}(0)x}L(0)$  is symmetric and positive definite. We denote it as following block matrix $\mathbf{M}$
\begin{align}
\mathbf{M}(x)\triangleq \left(\begin{array}{cc}
     \mathbf{M}_{11}(x) & \mathbf{M}_{12}(x)\\
     \mathbf{M}^T_{12}(x) & \mathbf{M}_{22}(x)
\end{array}\right)
\end{align}
according to  the block matrix $B_V(0)$.  Then by Cauchy-Schwarz inequality,  we get for all $\varepsilon>0$, that
\begin{align}\label{III}
\int^1_0 2V^T L^T(0)e^{-\mathbf{\Lambda}(0)x}L(0)B_V(0)V\,\mathrm{d}x
= &\int^1_0 2 v_1^T \mathbf{M}_{12}(x) \mathbf{S}(0)v_2 + 2v_2^T\mathbf{M}_{22}(x)\mathbf{S}(0)v_2\,\mathrm{d}x \nonumber\\
%2v_2^T\Big(S^T(0)\mathbf{M}_{22}(x)+\mathbf{M}_{22}(x)S(0)\Big)v_2\,\mathrm{d}x \nonumber\\
%\leq & C_{12}|v_2||v_1| + C_{22}||v_2||^2, \nonumber\\
\leq & \varepsilon||v_1||^2 + C_{\varepsilon}||v_2||^2,
\end{align}
where $C_{\varepsilon}>0$ is constant depending on $\varepsilon$.
Because $L^T(0)e^{-\mathbf{\Lambda}(0)x}\mathbf{\Lambda}^2(0)L(0)$ is   positive definite, there exists a constant $c_0>0$ such that
\begin{align}\label{IV}
-\int^1_0 V^T L^T(0)e^{-\mathbf{\Lambda}(0)x}\mathbf{\Lambda}^2(0)L(0)V \,\mathrm{d}x \leq -c_0 ||V||^2.
\end{align}

Thus, it follows from   \eqref{J0},\eqref{III}-\eqref{IV} and integrations by parts that
%\begin{align}\label{I}
%I\leq  -||v_2||^2  - V^T A_0(0)A(0)V|^1_0 + \int^1_0  O(V^3 + V^2V_t) \mathrm{d}x
%\end{align}
\begin{align}\label{IIe}
\mathbb{J}_0 \leq & (\varepsilon - c_0)||v_1||^2 + (C_{\varepsilon} - c_0) ||v_2||^2 +\Big[-V^T L^T(V) e^{-\mathbf{\Lambda}(0)x}\Lambda(V) L(V)V\Big]\Big|_0^1  \nonumber\\
& +\mathcal{O}\Big(\int^1_0 (|V|^3+|V|^2|V_t|) \,\mathrm{d}x; |V(t,\cdot)|_0\Big).
\end{align}

Let $\varepsilon = \frac{c_0}{2}$. Then combining \eqref{L0'2}, \eqref{I1} and \eqref{IIe} yields that
there exist positive constants $\alpha_0$,  $\beta_0$ independent of $V$, such that
\begin{align}\label{L01}
\mathbb{L}_0'(t) \leq  - \alpha_0 ||v_1||^2 + (\beta_0 - \alpha)||v_2||^2 +\mathbb{B}_0+ \mathcal{O}\Big(\int^1_0 (|V|^3+|V|^2|V_t|) \,\mathrm{d}x; |V(t,\cdot)|_0\Big),
\end{align}
where  the boundary term
\begin{align}\label{bc01}
\mathbb{B}_0 &= -V^T(t,x)\Big(\alpha A_0(V)A(V) + L^T(V)e^{-\mathbf{\Lambda}(0)x}\Lambda(V)L(V)\Big)V(t,x)\Big|^{1}_{0}.
\end{align}

It remains to  estimate $\mathbb{B}_0$. Using \eqref{matrixA2}, \eqref{bcvariable1} and the linear approximation, we have
\begin{align}\label{bc02}
\mathbb{B}_0 &= -\xi^T(t, x)\Big(\alpha (L^{-1}(0))^TA_0(0)L^{-1}(0)\mathbf{\Lambda}(0) +e^{-\mathbf{\Lambda}(0)x}\mathbf{\Lambda}(0)\Big)\xi(t, x)\Big|^{1}_{0} \nonumber \\
&\quad +\mathcal{O} \Big(|V(t,0)|^3+|V(t,1)|^3 ; |V(t,0)|+|V(t,1)| \Big).
\end{align}
Noting \eqref{matrixL} and \eqref{A0V}, \eqref{block-diagonal},  we easily obtain  that
\begin{align}\label{block-diagonal2}
(L^{-1}(0))^TA_0(0)L^{-1}(0) =
(\mathbf{L}^{-1}(0))^T\mathbf{A}_0(0)\mathbf{L}^{-1}(0) = \left(\begin{array}{cc}
     \mathbf{X}_1(0) & 0\\
     0 & \mathbf{X}_2(0)
\end{array}\right).
\end{align}
Substituting the boundary condition \eqref{bc}  and \eqref{block-diagonal2} into \eqref{bc02}, we thus get
\begin{align}\label{bc03}
\mathbb{B}_0 =
-\left(\begin{array}{c}
\xi_{+}(t, 1)\\
\xi_{-}(t, 0)
\end{array}\right)^T\mathbf{G}
\left(\begin{array}{c}
\xi_{+}(t, 1)\\
\xi_{-}(t, 0)
\end{array}\right)+\mathcal{O}\Big(|V(t,0)|^3+|V(t,1)|^3 ; |V(t,0)|+|V(t,1)| \Big)
\end{align}
where  $\mathbf{G}$ is a symmetric matrix defined as
\begin{align} \nonumber
\mathbf{G} \triangleq
& \,\alpha \left[\left(\begin{array}{cc}
      \mathbf{X}_2(0)\mathbf{\Lambda}_+(0) & 0\\
      0 &- \mathbf{X}_1(0)\mathbf{\Lambda}_-(0)
\end{array}\right)- \mathbf{K}^T\left(\begin{array}{cc}
     \mathbf{X}_2(0)\mathbf{\Lambda}_+(0) & 0\\
      0 & -\mathbf{X}_1(0)\mathbf{\Lambda}_-(0)
\end{array}\right)\mathbf{K}\right]  \\ \label{G}
 &+\left[\left(\begin{array}{cc}
     e^{-\mathbf{\Lambda}_+(0)}\mathbf{\Lambda}_+(0) & 0\\
      0 & -\mathbf{\Lambda}_-(0)
\end{array}\right)- \mathbf{K}^T\left(\begin{array}{cc}
     \mathbf{\Lambda}_+(0) & 0\\
      0 & -e^{-\mathbf{\Lambda}_-(0)}\mathbf{\Lambda}_-(0)
\end{array}\right)\mathbf{K}\right].
\end{align}
Using \eqref{bcvariable1}  and the boundary condition \eqref{bc}, we have
\begin{align} \label{Vxi1}
V(t,0)=L^{-1}(0)\xi(t,0)=L^{-1}(0)\left(\begin{array}{cc}
0 & \mathbf{I}_m \\
\mathbf{K}_{00}  &  \mathbf{K}_{01}
\end{array}\right)
\left(\begin{array}{c}
\xi_{+}(t, 1)\\
\xi_{-}(t, 0)
\end{array}\right)\\  \label{Vxi2}
V(t,1)=L^{-1}(0)\xi(t,1)=L^{-1}(0)
\left(\begin{array}{cc}
\mathbf{K}_{10} & \mathbf{K}_{11} \\
\mathbf{I}_{n-m}  &  0
\end{array}\right)\left(\begin{array}{c}
\xi_{+}(t, 1)\\
\xi_{-}(t, 0)
\end{array}\right)
\end{align}
Consequently, \eqref{bc03} becomes
\begin{align}
\mathbb{B}_0 =
-\left(\begin{array}{c}
\xi_{+}(t, 1)\\
\xi_{-}(t, 0)
\end{array}\right)^T
\mathbf{G}
\left(\begin{array}{c}
\xi_{+}(t, 1)\\
\xi_{-}(t, 0)
\end{array}\right)
+\mathcal{O}\Big((|V(t,0)|+|V(t,1)|)
\left|\left(\begin{array}{c}
\xi_{+}(t, 1)\\
\xi_{-}(t, 0)
\end{array}\right)\right|^2; |V(t,0)|+|V(t,1)|\Big).
\end{align}

Note that the assumptions \eqref{K1}-\eqref{K2} on the boundary feedback matrix $\mathbf{K}$ implies the symmetric matrix $\mathbf{G}$ is positive definite.
Therefore, there exist $\delta_{0}>0$ and $\gamma_0>0$ such that the boundary term $\mathbb{B}_0 \leq 0,$
and furthermore the estimate \eqref{el0} holds  if $|V(t,\cdot)|_0<\delta_{0}$. This concludes the proof of Lemma \ref{lem1}.

\subsection{Proof of Lemma 2.2}

By  \eqref{L1}, the time-derivative of $\mathbb{L}_1(t)$  can be expressed as
\begin{align}\label{L1t}
\mathbb{L}_1'(t) &= \int^1_0 2\Big(\alpha V_{t}^T A_0(V)V_{tt} + V^T_{t} L^T(V)e^{-\mathbf{\Lambda}(0)x}L(V) V_{tt}\Big)\,\mathrm{d}x
+ \mathcal{O}\Big(\int^1_0 |V_t|^3\,\mathrm{d}x; |V(t,\cdot)|_0\Big).
\end{align}
Substituting the term of $V_{tt}$ derived from \eqref{Vtt} into \eqref{L1t}, we have
\begin{align*}
\mathbb{L}_1'(t)
&= \int^1_0 2\alpha V_{t}^T A_0(V)\Big(B_V(V)V_{t}- (A'(V)V_t)V_x - A(V)V_{tx}\Big) \,\mathrm{d}x\\
&\quad + \int^1_02V^T_{t} L^T(V)e^{-\mathbf{\Lambda}(0)x}L(V) \Big(B_V(V)V_t - (A'(V)V_t)V_x - A(V)V_{tx}\Big)\,\mathrm{d}x\\
&\quad + \mathcal{O}\Big(\int^1_0 |V_t|^3\,\mathrm{d}x; |V(t,\cdot)|_0\Big)
\end{align*}
Thus, by using integrations by parts and some straightforward calculations, we get
\begin{align*}
\mathbb{L}_1'(t)
&= \int^1_0 \alpha V_t^T \Big(A_0(0)B_V(0) +B_V(0)^T A_0(0) \Big)V_t - V_{t}^T L^T(0)e^{-\mathbf{\Lambda}(0)x} \mathbf{\Lambda}^2(0) L(0)V_{t}\,\mathrm{d}x  \\
&\quad  +\int^1_0  2V_t^TL^T(0)e^{-\mathbf{\Lambda}(0)x}L(0)B_V(0)V_t\,\mathrm{d}x 
- V_t^T \Big(\alpha A_0(V)A(V) +  L^T(V) e^{-\mathbf{\Lambda}(0)x}\Lambda(V) L(V)\Big)V_{t}\Big|_0^1 \\
&\quad +\mathcal{O}\Big(\int^1_0 (|V_t|^3+|V||V_t|^2) \,\mathrm{d}x; |V(t,\cdot)|_0\Big).
\end{align*}
Similarly as the analysis of $\mathbb{L}_0'(t)$ in the proof of Lemma \ref{lem1}, we obtain that there exist positive constants $\alpha_1$ and $\beta_1$,
such that
\begin{align}
\mathbb{L}_1'(t)\leq  -\alpha_1||v_{1t}||^2 + (\beta_1 - \alpha)||v_{2t}||^2 + \mathbb{B}_1 + \mathcal{O}\Big(\int^1_0 (|V_t|^3+|V||V_t|^2) \,\mathrm{d}x; |V(t,\cdot)|_0\Big),
\end{align}
where the boundary term $\mathbb{B}_1$ is
\begin{align}\label{bc11}
\mathbb{B}_1 &= -\xi^T_t(t, x) \Big(\alpha (L^{-1}(0))^TA_0(0)L^{-1}(0)\mathbf{\Lambda}(0) +e^{-\mathbf{\Lambda}(0)x}\mathbf{\Lambda}(0)\Big)\xi_t(t, x)\Big|^{1}_{0}\nonumber\\
&\quad +\mathcal{O}\Big((|V(t,0)||V_t(t,0)|^2+|V(t,1)||V_t(t,1)|^2); |V(t,0)|+|V(t,1)| \Big).
\end{align}
Taking the time-derivative  \eqref{bc} and   \eqref{Vxi1}-\eqref{Vxi2}, we can easily express $(V_t(t,0), V_t(t,1))$  in terms of $(\xi_{+t}(t, 1), \xi_{-t}(t, 0))$, thus
the boundary term $\mathbb{B}_1$ can be rewritten as
\begin{align}
\mathbb{B}_1 &=
-\left(\begin{array}{c}
\xi_{+t}(t, 1)\\
\xi_{-t}(t, 0)
\end{array}\right)^T
\mathbf{G}
\left(\begin{array}{c}
\xi_{+t}(t, 1)\\
\xi_{-t}(t, 0)
\end{array}\right)\nonumber \\
&\quad +\mathcal{O}\Big((|V(t,0)|+|V(t,1)|)
\left|\left(\begin{array}{c}
\xi_{+t}(t, 1)\\
\xi_{-t}(t, 0)
\end{array}\right) \right|^2; |V(t,0)|+|V(t,1)| \Big).
\end{align}

Note again with  the assumptions \eqref{K1}-\eqref{K2} can imply  that $\mathbf{G}$ is positive definite.
Therefore, there exists $\delta_{1}>0$  and $\gamma_1>0$ such that $\mathbb{B}_1 \leq 0,$
and furthermore the estimate \eqref{el1} holds  if $|V(t,\cdot)|_0<\delta_{1}$. The finishes the proof of Lemma \ref{lem2}.

\subsection{Proof of Lemma 2.3}

Calculating the time-derivative of $\mathbb{L}_2(t)$ gives
\begin{align}
\mathbb{L}_2'(t) = \int^1_0 2\Big(\alpha V_{tt}^TA_0(V)V_{ttt} +  V^T_{tt} L^T(V)e^{-\mathbf{\Lambda}(0)x}L(V) V_{ttt}\Big)\,\mathrm{d}x
+ \mathcal{O}\Big(\int^1_0 |V_t||V_{tt}|^2\,\mathrm{d}x; |V(t,\cdot)|_0\Big).
\end{align}
Differentiating system \eqref{Vtt} with respect to $t$ and combining \eqref{systemV} and \eqref{Vtt}, we have,
\begin{align}\label{Vttt}
V_{ttt} + A(V)V_{ttx} = B_V(V)V_{tt} + (B_V(V))_{t}V_{t} -2(A'(V)V_t)V_{tx} -(A'(V)V_t)_t V_{x}.
\end{align}
Substituting the term of $V_{ttt}$ derived from \eqref{Vttt} into $\mathbb{L}_2'(t)$,
we do integration by parts and linear approximation, as in the proof of Lemma \ref{lem1} and \ref{lem2}, to deduce that
\begin{align}
\mathbb{L}_2'(t)
&= \int^1_0  \alpha V_{tt}^T \Big(A_0(0)B_V(0) + B_V(0)^TA_0(0)\Big) V_{tt} - V_{tt}^T L^T(0)e^{-\mathbf{\Lambda}(0)x} \mathbf{\Lambda}^2(0) L(0)V_{tt} \,\mathrm{d}x\nonumber\\
&\quad+  \int^1_0 2V_{tt}^T L^T(0)e^{-\mathbf{\Lambda}(0)x}L(0) B_V(0) V_{tt}\,\mathrm{d}x
 - V_{tt}^T\Big(\alpha A_0(V)A(V) + L^T(V) e^{-\mathbf{\Lambda}(0)x} \Lambda(V)L(V)\Big)V_{tt}  \Big|_0^1\nonumber\\
&\quad + \mathcal{O}\Big(\int^1_0 (|V||V_{tt}|^2
+ |V_t||V_{tt}|^2 + |V_t|^2 |V_{tt}|\,\mathrm{d}x; |V(t,\cdot)|_1\Big).
\end{align}
Thanks to the partially dissipative structure \eqref{matrixBV}-\eqref{A0B} and the fact that  $L^T(0)e^{-\mathbf{\Lambda}(0)x} \mathbf{\Lambda}^2(0) L(0)$ is positive definite,
there exist positive constants $\alpha_2$ and $\beta_2$  independent of $V$ such that
\begin{align} \nonumber
\mathbb{L}_2'(t)& \leq  -\alpha_2||v_{1tt}||^2 + (\beta_2- \alpha)||v_{2tt}||^2  + \mathbb{B}_2
\\
&\quad +  \mathcal{O}\Big(\int^1_0 (|V||V_{tt}|^2
+ |V_t||V_{tt}|^2 + |V_t|^2 |V_{tt}|\,\mathrm{d}x; |V(t,\cdot)|_1\Big)
  \end{align}
 where  $\mathbb{B}_2$  denotes the boundary term derived from integration by parts.
Taking the second time-derivative of \eqref{bc} and   \eqref{Vxi1}-\eqref{Vxi2}, we can rewrite the boundary term $\mathbb{B}_2$ as
\begin{align}\label{bc21}
\mathbb{B}_2 &=
-\left(\begin{array}{c}
\xi_{+tt}(t, 1)\\
\xi_{-tt}(t, 0)
\end{array}\right)^T
\mathbf{G}
\left(\begin{array}{c}
\xi_{+tt}(t, 1)\\
\xi_{-tt}(t, 0)
\end{array}\right)\nonumber \\
&\quad +\mathcal{O}\Big((|V(t,0)|+|V(t,1)|)
\left|\left(\begin{array}{c}
\xi_{+tt}(t, 1)\\
\xi_{-tt}(t, 0)
\end{array}\right)\right|^2; |V(t,0)|+|V(t,1)|\Big).
\end{align}
Note again that $\mathbf{G}$ is positive definite if  the assumptions \eqref{K1}-\eqref{K2} hold.
Therefore, there exists $\delta_{2}>0$  and $\gamma_2>0$ such that $\mathbb{B}_2 \leq 0$
and furthermore the estimate \eqref{el2} holds  if $|V(t,\cdot)|_1<\delta_{2}$. This ends the proof of Lemma \ref{lem3}.

\section{Application to  Saint-Venant-Exner equations}

We now consider the Saint-Venant-Exner equations for a moving bathymetry on a
sloping channel with a rectangular cross-section:
\begin{align}\label{SVE}
\left.\begin{array}{lll}
&H_t + VH_x + HV_x = 0,\\
&B_t + aV^2V_x = 0,\\
&V_t + VV_x + gH_x + gB_x = gS_b - C_f\frac{V^2}{H}.
\end{array}\right.
\end{align}
Here $H = H(t, x)$ is the water depth, $B=B(t,x)$ is the elevation of the sediment bed and $V = V(t, x)>0$ is the average velocity of water. Moreover, $g$ is the gravity constant, constant $S_b$ is the bottom slope of the channel, constant $C_f$ means the friction coefficient and constant $a$ is a parameter that includes porosity and viscosity effects on the sediment dynamics (see \cite{Hudson2003}).

Let   $(H_*,B_*, V_*)^T$ with $H_*>0, B_*>0$ and $ V_*>0$ be a  constant equilibrium of  \eqref{SVE}, i.e.,
\begin{align}\label{steady2}
gS_bH_* = C_f(V_*)^2 > 0.
\end{align}
Let $U \triangleq  (h, b, v)^T$ be  the deviations of the states:
\begin{align}
h(t,x)=H(t,x)-H_*,\quad
b(t,x)=B(t,x)-B_*,\quad
v(t,x)=V(t,x)-V_*.
\end{align}
Then the Saint-Venant-Exner equations \eqref{SVE} can be rewritten in the form of \eqref{system}
with
\begin{align}\label{SVEA}
\mathbf{A}(U) = \left(\begin{array}{ccc}
     v+V_* & 0 & h+H_* \\
     0 & 0 & a(v+V_*)^2\\
     g & g & v+V_*
\end{array}\right),
\quad \mathbf{Q}(U)
= \left(\begin{array}{c}
    0\\
    0\\
gS_b - C_f\frac{(v+V_*)^2}{h+H_*}
\end{array}\right).
\end{align}

In order to apply Theorem \ref{thm}, we will verify that the hyperbolic system satisfies the \emph{partially dissipative structure}.  For simplicity,
we only show that  \eqref{eigenvalues}, \eqref{hyperbolic}, \eqref{structural1} and \eqref{structural2}  are satisfied at the equilibrium $U=0$.

First, the matrix $\mathbf{A}(0)$  has three eigenvalues $\lambda_i$ $(i=1,2,3)$ satisfying
\begin{align}\label{SVEeigenvalue}
\lambda^3-2V_*\lambda^2+(V_*^{2}-gaV_*^{2}-gH_*)\lambda+gaV_*^{3}=0.
\end{align}
Thus we get the following relations
\begin{align}\label{SVEeigenvalues2}
\left.\begin{array}{lll}
\lambda_1\lambda_2\lambda_3=-gaV_*^{3},\\
\lambda_1+\lambda_2+\lambda_3=2V_*,\\
\lambda_1\lambda_2+\lambda_1\lambda_3+\lambda_2\lambda_3=V_*^{2}-gaV_*^{2}-gH_*.\end{array}\right.
\end{align}
Here $\lambda_1$, $\lambda_3$ are the characteristic velocities of the water flow and $\lambda_2$ is
the characteristic velocity of the sediment motion. Therefore, \eqref{steady2} and \eqref{SVEeigenvalues2} yield that matrix $\mathbf{A}(0)$ has no vanishing eigenvalues with
\begin{align} \label{lambda123}
\lambda_1<0<\lambda_2\leq\lambda_3.
\end{align}

Due to the fact that the sediment motion is much slower than the water flow, we
can make the following reasonable assumptions that $\lambda_2$ is so small that
\begin{align}\label{lambda2small}
\lambda_2<-\lambda_1\quad \text{and}\quad 0<\lambda_2<\frac{3}{2}V_*,
\end{align}
which leads to the following relations
\begin{align}\label{elambda123}
\lambda_1-\frac{3}{2}V_*<0, \quad \lambda_2-\frac{3}{2}V_*<0,\quad \lambda_3-\frac{3}{2}V_*=\frac{1}{2}V_*-\lambda_1-\lambda_2>0.
\end{align}

Next, we choose an invertible matrix $\mathbf{P}(0)$ and a symmetric positive definite matrix $\mathbf{A}_0(0)$, such that \eqref{hyperbolic},  \eqref{structural1} and \eqref{structural2} are satisfied. Inspired by \cite{Yong2016}, we take
\begin{align}
\mathbf{P}(0) = \left(\begin{array}{ccc}
     1 & 0 & 0\\
     0 & 1 & 0\\
     -\frac{V_*}{2H_*} & 0 & 1
\end{array}\right)
\text{\quad and \quad}
\mathbf{A}_0(0)=\left(\begin{array}{ccc}
     \frac{4gH_*+2agV_*^{2}}{4H_*^{2}} & -\frac{g}{2H_*} & -\frac{V_*}{2H_*}\\
     -\frac{g}{2H_*} & \frac{3g}{2aV_*^{2}} & 0\\
     -\frac{V_*}{2H_*} & 0 & 1
\end{array}\right).
\end{align}

According to \eqref{SVEeigenvalues2} and \eqref{elambda123}, it is easy to verify that
\begin{align}\label{i31}
\det(\mathbf{A}_0(0))=\frac{3g^2}{2aH_*V_*^{2}}+\frac{g^2}{2H_*^{2}}-\frac{3g}{8a^2H_*^{2}}
=\frac{g}{aH_*^{2}V_*^{3}}\prod_{i=1}^3\Big(\lambda_i-\frac{3}{2}V_*\Big)>0.
\end{align}
Therefore, $\mathbf{A}_0(0)$ is symmetric positive definite.
It  follows also that $\mathbf{A}_0(0)\mathbf{A}(0)=\mathbf{A}^T(0)\mathbf{A}_0(0)$.

Note that
\begin{align}\label{SVEQU}
\mathbf{Q}_U(0) = \left(\begin{array}{ccc}
     0 & 0 & 0\\
     0 & 0 & 0\\
     C_f\frac{{V_*}^2}{{H_*}^2} & 0 & -2C_f\frac{V_*}{H_*}
\end{array}\right).
\end{align}
Direct computations give that
\begin{align}
&\mathbf{P}(0)\mathbf{Q}_U(0)\mathbf{P}^{-1}(0) = \left(\begin{array}{ccc}
     0 & 0 & 0\\
     0 & 0 & 0\\
     0 & 0 & -C_f\frac{V_*}{H_*}
\end{array}\right)
\\
&(\mathbf{P}^{-1}(0))^T\Big(\mathbf{A}_0(0)\mathbf{Q}_U(0) + \mathbf{Q}_U^T(0)\mathbf{A}_0(0)\Big)\mathbf{P}^{-1}(0)=
\left(\begin{array}{ccc}
     0 & 0 & 0\\
     0 & 0& 0\\
     0 & 0& -4C_f\frac{V_*}{H_*}
\end{array}\right).
\end{align}
i.e., the \emph{partially dissipative structure}  indeeds holds.

Let $\mathbf{L}_i(0)$  be the left eigenvectors corresponding to $\lambda_i$ $(i=1,2,3)$, then we have
\begin{align}
\mathbf{L}(0)=
\left(\begin{array}{c}
\mathbf{L}_1(0)\\
\mathbf{L}_2(0)\\
\mathbf{L}_3(0)
\end{array}\right)
=
\left(\begin{array}{ccc}
\frac{g}{\lambda_1-V_*} & \frac{g}{\lambda_1} &1\\
\frac{g}{\lambda_2-V_*} & \frac{g}{\lambda_2} &1\\
\frac{g}{\lambda_3-V_*} & \frac{g}{\lambda_3} &1
\end{array}\right)
\end{align}
and further
\begin{align}
\mathbf{L}^{-1}(0)=\left(\begin{array}{ccc}
\frac{\lambda_1H_*}{(\lambda_1-\lambda_2)(\lambda_1-\lambda_3)} &
\frac{\lambda_2H_*}{(\lambda_2-\lambda_3)(\lambda_2-\lambda_1)} &
\frac{\lambda_3H_*}{(\lambda_3-\lambda_1)(\lambda_3-\lambda_2)}   \vspace{2mm} \\
\frac{aV_*^{2}(\lambda_1-V_*)}{(\lambda_1-\lambda_2)(\lambda_1-\lambda_3)} & \frac{aV_*^{2}(\lambda_2-V_*)}{(\lambda_2-\lambda_3)(\lambda_2-\lambda_1)} &
\frac{aV_*^{2}(\lambda_3-V_*)}{(\lambda_3-\lambda_1)(\lambda_3-\lambda_2)}  \vspace{2mm} \\
\frac{\lambda_1(\lambda_1-V_*)}{(\lambda_1-\lambda_2)(\lambda_1-\lambda_3)} & \frac{\lambda_2(\lambda_2-V_*)}{(\lambda_1-\lambda_2)(\lambda_1-\lambda_3)} & \frac{\lambda_3(\lambda_3-V_*)}{(\lambda_3-\lambda_1)(\lambda_3-\lambda_2)}
\end{array}\right).
\end{align}
Consequently, we have
\begin{align}
(\mathbf{L}^{-1}(0))^T\mathbf{A}_0(0)\mathbf{L}^{-1}(0)
=\left(\begin{array}{ccc}
    \mathbf{X}_{11} & 0   & 0\\
     0 & \mathbf{X}_{22} & 0\\
    0 & 0 & \mathbf{X}_{33}
\end{array}\right)
\end{align}
where
\begin{align}
 \mathbf{X}_{11} = \frac{\lambda_{1}(\lambda_1-\frac{3}{2}V^{*})}{(\lambda_1-\lambda_2)(\lambda_1-\lambda_3)},
 \quad  \mathbf{X}_{22}= \frac{\lambda_{2}(\lambda_2-\frac{3}{2}V^{*})}{(\lambda_2-\lambda_1)(\lambda_2-\lambda_3)}
 \quad  \mathbf{X}_{33}= \frac{\lambda_{3}(\lambda_3-\frac{3}{2}V^{*})}{(\lambda_3-\lambda_1)(\lambda_3-\lambda_2)}.
\end{align}

For the Saint-Venant-Exner equations \eqref{SVE}, we introduce the following physical boundary conditions
\begin{align}\label{SVEbc}
\left\{\begin{array}{lll}
b(t,0)=0,\\
v(t,0)=-k_1h(t,0),\\
v(t,1)=-k_2(h(t,1)+b(t,1)).
\end{array}\right.
\end{align}
where  $k_1$ and $k_2$ are two feedback parameters.

Then the boundary conditions \eqref{SVEbc} can be rewritten in the form of its linearized diagonal variable $\xi$ defined by \eqref{bcvariable} as
\begin{align}\label{xibc}
\left(\begin{array}{ccc}
\xi_2(t,0)\\
\xi_3(t,0)\\
\xi_1(t,1) \end{array}\right) =
\mathbf{K}
\left(\begin{array}{ccc}
\xi_2(t,1)\\
\xi_3(t,1)\\
\xi_1(t,0) \end{array}\right)
\text{\quad with \quad}
\mathbf{K}=
\left(\begin{array}{ccc}
0&0&\pi_2(k_1)\\
0&0&\pi_3(k_1)\\
\chi_2(k_2)&\chi_3(k_2)&0
\end{array}\right),
\end{align}
where $\pi_j$ and $\chi_j$ $(j=2,3)$ are the following quantities depending on $k_1$ and $k_2$:
\begin{align}
&\pi_j(k_1)=\frac{\lambda_1-V_*}{\lambda_j-V_*} \cdot \frac{g-k_1(\lambda_j-V_*)}{g-k_1(\lambda_1-V_*)}\quad (j=2,3),\\
&\chi_2(k_2)=\frac{\lambda_2(\lambda_3-\lambda_1)(\lambda_2-V_*)}{\lambda_1(\lambda_3-\lambda_2)(\lambda_1-V_*)}
  \cdot \frac{g+k_2(\lambda_2-V_*)}{g+k_2(\lambda_1-V_*)},\\
&\chi_3(k_2)=\frac{\lambda_3(\lambda_1-\lambda_2)(\lambda_3-V_*)}{\lambda_1(\lambda_3-\lambda_2)(\lambda_1-V_*)}
  \cdot \frac{g+k_2(\lambda_3-V_*)}{g+k_2(\lambda_1-V_*)}.
\end{align}

Consequently, the dissipative boundary conditions \eqref{K1} and \eqref{K2} yield
\begin{align}\label{SVEk11}
&\pi_2^2(k_1)\frac{\mathbf{X}_{22}}{\mathbf{X}_{11}}\frac{\lambda_2}{|\lambda_1|} + \pi_3^2(k_1)\frac{\mathbf{X}_{33}}{\mathbf{X}_{11}}\frac{\lambda_3}{|\lambda_1|}\leq 1,\\ \label{SVEk21}
&\chi_2^2(k_2)\frac{\mathbf{X}_{11}}{\mathbf{X}_{22}}\frac{|\lambda_1|}{\lambda_2}+ \chi_3^2(k_2)\frac{\mathbf{X}_{11}}{\mathbf{X}_{33}}\frac{|\lambda_1|}{\lambda_3}\leq 1
\end{align}
and
\begin{align}\label{SVEk12}
& \pi_2^2(k_1)\frac{\lambda_2}{|\lambda_1|} + \pi_3^2(k_1)\frac{\lambda_3}{|\lambda_1|}\leq 1,\\ \label{SVEk22}
&\chi_2^2(k_2)e^{\lambda_2-\lambda_1}\frac{|\lambda_1|}{\lambda_2}+ \chi_3^2(k_2)e^{\lambda_3-\lambda_1}\frac{|\lambda_1|}{\lambda_3}\leq 1.
\end{align}

Let
\begin{align}\label{beta23}
\beta_j = \max\Big\{\frac{\mathbf{X}_{jj}}{\mathbf{X}_{11}},1\Big\}, \
\eta_j =\max\Big\{\frac{\mathbf{X}_{11}}{\mathbf{X}_{jj}},e^{\lambda_j-\lambda_1}\Big\},
\quad (j=2,3).
\end{align}
Then,  \eqref{SVEk11}, \eqref{SVEk21}, \eqref{SVEk12} and \eqref{SVEk22} are satisfied if
\begin{align}\label{SVEk13}
\pi_2^2(k_1)\beta_2\frac{\lambda_2}{|\lambda_1|} + \pi_3^2(k_1)\beta_3\frac{\lambda_3}{|\lambda_1|}\leq 1,\\ \label{SVEk23}
\chi_2^2(k_2)\eta_2\frac{|\lambda_1|}{\lambda_2} + \chi_3^2(k_2)\eta_3\frac{|\lambda_1|}{\lambda_3}\leq 1.
\end{align}

Finally we conclude by Theorem \ref{thm} that
\begin{thm}\label{thm app}
If the boundary feedback parameters $k_1$ and $k_2$ satisfy \eqref{SVEk13} and \eqref{SVEk23}, respectively,
then the constant equilibrium $(H_*,B_*, V_*)^T$ of the Saint-Venant-Exner system \eqref{SVE}, \eqref{SVEbc} is locally exponentially
stable for the $H^2$-norm.
\end{thm}

\section*{Acknowledgements}

The authors were partially supported by the National Science Foundation of China (No. 11971119) and the Fundamental Research Funds for the Central Universities (No. 2232020D-41).

\bibliographystyle{plain}
\bibliography{mybib}

\begin{thebibliography}{10}

\bibitem{BC11}
Georges Bastin and Jean-Michel Coron.
\newblock On boundary feedback stabilization of non-uniform linear {$2\times2$}
  hyperbolic systems over a bounded interval.
\newblock {\em Systems Control Lett.}, 60(11):900--906, 2011.

\bibitem{Coron2016book}
Georges Bastin and Jean-Michel Coron.
\newblock {\em Stability and Boundary Stabilization of 1-{D} Hyperbolic
  Systems}, volume~88 of {\em Progress in Nonlinear Differential Equations and
  their Applications}.
\newblock Birkh\"{a}user/Springer, [Cham], 2016.
\newblock Subseries in Control.

\bibitem{CoronBook07}
Jean-Michel Coron.
\newblock {\em Control and nonlinearity}, volume 136 of {\em Mathematical
  Surveys and Monographs}.
\newblock American Mathematical Society, Providence, RI, 2007.

\bibitem{Coron2008dissipative}
Jean-Michel Coron, Georges Bastin, and Brigitte d'Andr\'{e}a Novel.
\newblock Dissipative boundary conditions for one-dimensional nonlinear
  hyperbolic systems.
\newblock {\em SIAM J. Control Optim.}, 47(3):1460--1498, 2008.

\bibitem{CN15}
Jean-Michel Coron and Hoai-Minh Nguyen.
\newblock Dissipative boundary conditions for nonlinear 1-{D} hyperbolic
  systems: sharp conditions through an approach via time-delay systems.
\newblock {\em SIAM J. Math. Anal.}, 47(3):2220--2240, 2015.

\bibitem{Coron2013backstepping}
Jean-Michel Coron, Rafael Vazquez, Miroslav Krstic, and Georges Bastin.
\newblock Local exponential {$H^2$} stabilization of a {$2\times 2$}
  quasilinear hyperbolic system using backstepping.
\newblock {\em SIAM J. Control Optim.}, 51(3):2005--2035, 2013.

\bibitem{Dick2011strict}
Markus Dick, Martin Gugat, and G\"{u}nter Leugering.
\newblock A strict {$H^1$}-{L}yapunov function and feedback stabilization for
  the isothermal {E}uler equations with friction.
\newblock {\em Numer. Algebra Control Optim.}, 1(2):225--244, 2011.

\bibitem{Greenberg1984effect}
J.~M. Greenberg and Tatsien Li.
\newblock The effect of boundary damping for the quasilinear wave equation.
\newblock {\em J. Differential Equations}, 52(1):66--75, 1984.

\bibitem{Gugat2017}
Martin Gugat, G\"{u}nter Leugering, and Ke~Wang.
\newblock Neumann boundary feedback stabilization for a nonlinear wave
  equation: a strict {$H^2$}-{L}yapunov function.
\newblock {\em Math. Control Relat. Fields}, 7(3):419--448, 2017.

\bibitem{Hayat2018}
Amaury Hayat.
\newblock Boundary stability of 1-{D} nonlinear inhomogeneous hyperbolic
  systems for the {$C^1$} norm.
\newblock {\em SIAM J. Control Optim.}, 57(6):3603--3638, 2019.

\bibitem{Yong2016}
Michael Herty and Wen-An Yong.
\newblock Feedback boundary control of linear hyperbolic systems with
  relaxation.
\newblock {\em Automatica J. IFAC}, 69:12--17, 2016.

\bibitem{HRDK19}
Long Hu, Rafael Vazquez, Florent Di~Meglio, and Miroslav Krstic.
\newblock Boundary exponential stabilization of 1-dimensional inhomogeneous
  quasi-linear hyperbolic systems.
\newblock {\em SIAM J. Control Optim.}, 57(2):963--998, 2019.

\bibitem{Hudson2003}
Justin Hudson and Peter~K. Sweby.
\newblock Formulations for numerically approximating hyperbolic systems
  governing sediment transport.
\newblock {\em J. Sci. Comput.}, 19(1-3):225--252, 2003.
\newblock Special issue in honor of the sixtieth birthday of Stanley Osher.

\bibitem{Krstic2008book}
Miroslav Krstic and Andrey Smyshlyaev.
\newblock {\em Boundary control of {PDE}s}, volume~16 of {\em Advances in
  Design and Control}.
\newblock Society for Industrial and Applied Mathematics (SIAM), Philadelphia,
  PA, 2008.
\newblock A course on backstepping designs.

\bibitem{Li1994book}
Tatsien Li.
\newblock {\em Global Classical Solutions for Quasilinear Hyperbolic Systems},
  volume~32 of {\em RAM: Research in Applied Mathematics}.
\newblock Masson, Paris; John Wiley \& Sons, Ltd., Chichester, 1994.

\bibitem{Xu2002}
Cheng-Zhong Xu and Gauthier Sallet.
\newblock Exponential stability and transfer functions of processes governed by
  symmetric hyperbolic systems.
\newblock {\em ESAIM Control Optim. Calc. Var.}, 7:421--442, 2002.

\bibitem{Yong1999}
Wen-An Yong.
\newblock Singular perturbations of first-order hyperbolic systems with stiff
  source terms.
\newblock {\em J. Differential Equations}, 155(1):89--132, 1999.

\bibitem{Yong2019}
Wen-An Yong.
\newblock Boundary stabilization of hyperbolic balance laws with characteristic
  boundaries.
\newblock {\em Automatica J. IFAC}, 101:252--257, 2019.

\end{thebibliography}

\end{CJK*}

\end{document}